\documentclass[leqno,a4paper,11pt]{article}
\usepackage{amsmath}
\usepackage{amsfonts}
\usepackage{amssymb}
\newtheorem{theo}{Theorem}[section]

\newtheorem{cor}[theo]{Corollary}
\newtheorem{lem}[theo]{Lemma}
\newtheorem{remark}[theo]{Remark}

\newtheorem{Bestimmung}[theo]{Definition}%[section]

\newtheorem{example}[theo]{Example}
\newtheorem{examples}[theo]{Examples}

\newtheorem{counterexample}[theo]{Counterexample}

\newtheorem{remdefi}[theo]{Remark and Definition}

\newcommand\titre[1]{{\bf(#1)} } %% titre de theoreme

 \newcommand\preuv[1]{\par\noindent{\bf #1.}\ } %{\par{\sc #1.}\  }
 \newcommand\preuve{\preuv{Proof}}
 \newcommand\preuvof[1]{\preuv{Proof of {#1}}}

 \newcommand\carre{\fbox{\rule{0em}{.3em}\rule{.3em}{0em}} \qquad}
 \newcommand\fin{{\rule{.1cm}{0cm}}\hfill$\carre\qquad$\medskip\par}

\newcommand\C{\mathbb C}%
\newcommand\N{\mathbb N}%
\newcommand\R{\mathbb R}%
\newcommand\X{\mathbb X}%

\newcommand\abs[1]{\left\vert{#1}\right\vert}
\newcommand\accol[1]{\left\{{#1}\right\}}
\newcommand\CCO[1]{\left({#1}\right)}

\newcommand\norm[1]{\left\Vert{#1}\right\Vert}

\newcommand\indicatrix{\,\rlap{{\rm 1}}\kern.22em \mbox{\rm l}}
\newcommand\un[1]{\indicatrix_{#1}} %% fonction indicatrice de #1

\newcommand\restr[1]{ { \rule[-5pt]{0.4pt}{10pt}}_{#1} } % restriction to a
\newcommand\ra{\rightarrow}
\newcommand\mt{\mapsto}

 \newcommand\tq{;\,}
 \newcommand\foreach[1]{{\forall{#1}}\hspace{.8em}}
 \newcommand\thereis[1]{{\exists{#1}}\hspace{.8em}}

\newcommand\graphe[1]{\mathop{\text{\rm gph}}\CCO{#1}}     % graphe de #1

\newcommand\point{{\displaystyle{\text .}}} % pour obtenir un gros ``.'' dans
\newcommand\lsc{l.s.c.~}

\newcommand\projtheo{Projection Theorem}

\newcommand\proj[1]%{\mathop{\text{\rm proj}_{#1}}} %% projection
 {\mathop{\text{\rm proj}_{#1}}}

\newcommand\esp{\mathbb{S}}            %% THE topological space
\newcommand\espp{\esp_0}  % a subspace
\newcommand\esptopol{\mathfrak{T}}     %% some topological space
\newcommand\pol{\mathbb{T}} % a Polish space
\newcommand\bor[1]{ \mathcal{B}_{#1} }     %% Borel $\sigma$--algebra of #1
\newcommand\elt{x}                      %% an element of \esp
\newcommand\eltt{y}                      %% an other element of \esp
\newcommand\Cfcts{\mathop{\text{\rm C}}\nolimits} % symbole generique
\newcommand\Cb[1]{\Cfcts_b\CCO{#1}} % fcts continues bornees

\newcommand\dense{{D}} %% un sous-espace dense

\newcommand\compacta{{\mathcal K}} % ens des parties compactes de \esp

\newcommand\esprob{\Omega}                       %% espace probabilisable
\newcommand\tribu{{\mathcal A}}                  %% tribu
\newcommand\pr{ \mathop{\text{\rm P}}\nolimits }          %% probabilite
\newcommand\oap{(\esprob, \tribu, \pr)}        %% espace probabilise
\newcommand\elprob{\omega}

\newcommand\sstribu{{\mathcal M}} %% une sous-tribu

\newcommand\espmesgen{{\mathbb M}}
\newcommand\tribugen{{\mathcal M}}

\newcommand\laws[2][]{{\mathcal P_{#1}(#2)}}

\newcommand\complet[2][]{{#2}^*_{#1}}          %% complete de #1 (tribu ou proba)

\newcommand\espoir{\mathop{\kern 0pt{E}}\nolimits} %%symbole generique
\newcommand\expect[1]{{\espoir}\CCO{#1}}%%esperance mat.(cf aussi \condexp)
\newcommand\condexpect[2]{ {\espoir}\CCO{#1|{#2}} } %% esperance cond % la
\renewcommand\L{\text{\rm L}} %% espaces\L^p

\newcommand\cout{c}
\newcommand\DistrJointes{{D}} %% ens. des couplages de #1 et #2 
\newcommand\KR{{\Delta}_{{\text{\rm KR}}}^{(\cout)}  }
\newcommand\paramKR{\underline{\Delta}_{{\text{\rm KR}}}^{(\cout)}  }
\newcommand\paramDistrJointes{\underline{D}} %% ens. des couplages de #1 et #2

\newcommand\lip[2]{\mathop{ \text{\rm Lip}^{(#1)}_{\,#2} }} %fcts lipschitziennes 
\newcommand\Lipc{\lip{\cout}{\esp}}

\newcommand\lipdist[1]{\Delta_{\text{\rm L}}^{(#1)}}
\newcommand\lipdistc{\lipdist{\cout}}
\newcommand\paramLip[1]{\mathop{\underline{\text{\rm Lip}}^{(#1)}}} %fcts lipschitziennes 
\newcommand\paramLipc{\paramLip{\cout}}
\newcommand\Lipcm{\lip{\cout}{\esp, \cal{M}}}
\newcommand\paramlipdist[1]{\underline{\Delta}_{\text{\rm L}}^{(#1)}}
\newcommand\paramlipdistc{\paramlipdist{\cout}}
\newcommand\niveaubas{r} %%utilise dans lem \ref{KRmeas} et dans th \ref{KR}

\newcommand\youngs{{\mathcal Y}}
\newcommand\youngsint[1]{\youngs^{\,#1,1}} % mesures de Young
\author{J\'er\^ome {\sc Dedecker}\footnote{Laboratoire de Statistique Th\'eorique et Appliqu\'ee, Universit\'e
Paris 6,  175 rue du Chevaleret, 75013 Paris, France.
E-mail: dedecker@ccr.jussieu.fr}, 
        Cl\'ementine {\sc Prieur}\footnote{Laboratoire
de Statistique et Probabilit\'es, Universit\'e Paul Sabatier, 118 route de
Narbonne, 31062 Toulouse cedex 4, France. E-mail: prieur@cict.fr}, 
       Paul {\sc Raynaud de Fitte}\footnote{Laboratoire Rapha\"el Salem, UMR
   CNRS 6085, UFR Sciences , Universit\'e de Rouen, 76821 
   Mont Saint Aigan Cedex, France. E-Mail: prf@univ-rouen.fr}}
\title
{Parametrized Kantorovich-Rubin\v stein theorem  and application to the coupling of random variables \\
\medskip
Th\'eor\`eme de Kantorovich-Rubin\v stein avec param\`etre et application au couplage 
des variables al\'eatoires}
\date{}
\begin{document}
\maketitle

\begin{abstract}
We prove a version for random measures of the celebrated Kantorovich-Rubin\v stein duality
theorem and we give an application to the coupling
of random variables which extends and unifies known results. 
\begin{center}
{\bf R\'esum\'e} 
\end{center}

Nous d\'emontrons une version du  th\'eor\`eme de dualit\'e de  Kantorovich-Rubin\v stein pour les mesures
al\'eatoires, et nous donnons une application  au couplage des variables al\'eatoires qui
\'etend et unifie les r\'esulats ant\'erieurs.

\end{abstract}

%%%%%%%%%%%%%%%%%%%%%%%
\section{Introduction and notations}

Let  $\mu$ and $\nu$ be two probability measures
on a Polish space $({\mathbb S}, d)$. In 1970
Dobru{\v{s}}in \cite[page 472]{dobrushin70var} proved that there exists a probability measure 
$\lambda$ on ${\mathbb S } \times {\mathbb S}$  with margins $\mu$ and $\nu$, such that
\begin{equation}\label{dobr}
    \lambda(\{x\neq y, (x,y) \in {\mathbb S} \times {\mathbb S} \})=\frac{1}{2}\| \mu-\nu \|_v \, ,
\end{equation}
where $\|\cdot\|_v$ is the variation norm. More precisely, Dobru{\v{s}}in gave an explicit solution to (\ref{dobr}) defined by
\begin{equation}\label{solution}
  \lambda(A \times B)= (\mu-  \pi_-)(A \cap B) + \frac{\pi_-(A) \pi_+(B)}{\pi_+({\mathbb S})} \quad \text{for $A$, $B$ in 
${\mathcal B}_{\mathbb S}$}\, ,
\end{equation}
where $\mu-\nu= \pi_+-\pi_-$ is the Hahn decomposition of $\pi= \mu-\nu$.

%More precisely he proved that if $\lambda^+- \lambda^-$ is the Jordan
%decomposition of $P-Q$, then the measure 
%\begin{equation}\label{ex1}
%  \mu(A\times B)=P(A \cap B)- \lambda^+(A \cap B) + \frac{ \lambda^-(A) \lambda^+(B)}{\lambda^+({\cal X})}
%\end{equation}
%is solution of (\ref{dobr}).
Starting from (\ref{solution}) (see \cite[Proposition 4.2.1]{Berbee79book}), Berbee 
obtained the following coupling result (\cite[Corollary 4.2.5]{Berbee79book}):
let $(\Omega, {\cal A}, {\mathbb P})$ be a probability space, let ${\cal M}$ be a 
$\sigma$-algebra of ${\cal A}$, and let $X$ be a random variable
with values in ${\mathbb S}$. Denote by $\pr_X$ the distribution of $X$ and
by $\pr_{X|{\mathcal M}}$ a regular conditional distribution of $X$ given ${\mathcal M}$.
If $\Omega$ is rich enough, there exists $X^*$ distributed as $X$ 
and independent of ${\cal M}$ such that
\begin{equation}\label{ber}
   \pr({X\neq X^*})=\frac{1}{2}E( \|{\pr}_{X|{\cal M}} -{\pr}_X \|_v)  \, .
\end{equation}
To prove  (\ref{ber}), Berbee built  a couple $(X, X^*)$ whose conditional 
distribution given ${\mathcal M}$ is the random probability $\lambda_\omega$ defined by
(\ref{solution}), with  random 
margins $\mu={\pr}_{X|{\cal M}}$ and $\nu={\pr}_{X}$.

It is by now well known that Dobru{\v{s}}in's result (\ref{dobr}) is a particular case 
of the Kantorovich-Rubin\v stein duality theorem (which we recall  at the beginning of Section 2) 
applied to the discrete metric $c(x,y)=\indicatrix_{x \neq y}$ (see \cite[page 93]{rachev-ruschendorf98bookI}). 
Starting from this simple remark,
Berbee's proof can be described as follows:  one can find
a couple 
$(X, X^*)$ whose conditional distribution with respect to ${\mathcal M}$ solves
the duality problem with cost function  $c(x,y)=\indicatrix_{x \neq y}$ and random margins
$\mu=\pr_{X|{\mathcal M}}$ and  $\nu= \pr_X$. 

A reasonable question is then: for what class of cost functions can we obtain the 
same kind of coupling than Berbee's? Or, equivalently, given two random
probabilities $\mu_\omega$ and $\nu_\omega$ on a Polish space $({\mathbb S},d)$, for what class of cost functions
is there a random probability on 
${\mathbb S} \times {\mathbb S}$ solution to the duality problem with  margins $(\mu_\omega, \nu_\omega)$?
%The Kantorovitch-Rubin\v stein duality Theorem has a
%rich history and many applications. It first appeared in the 1958 work of
%Kantorovich and Rubin\v stein on the mass transport. A general form of this
%theorem is recalled at the beginning of Section 2. One can wonder what
%happens when considering random probabilities in this theorem? That is the purpose of this paper.
%More precisely: given two random
%probabilties $\mu_\omega$ and $\nu_\omega$ on a Polish space $({\mathbb S},d)$, for what class of cost functions
%is there a random probability on 
%${\mathbb S} \times {\mathbb S}$ solution to the duality problem with  margins $(\mu_\omega, \nu_\omega)$?
In 2004, this question has been partially answered in two independent works. In 
Proposition 1.2 of \cite{dedecker-prieur04coupl}
the authors prove the existence of such a random probability for 
the cost function $c=d$.
From the proof of Theorem 3.4.1 in \cite{cc-prf-valadier04book}, we see that this result holds in fact
for any  distance $c$ which is continuous with respect to $d$. To summarize,  we know that 
there exists a random probability solution to the duality problem with cost function $c$ 
and  given random
margins in two  distinct
situations: on  one hand $c$ is the discrete metric, on the other hand $c$ is
any continuous distance  with respect to $d$.
A general result containing both situations as particular cases would be more 
satisfactory.

The main result of this paper (point 1 of Theorem 2.1) 
asserts that there exists a random probability on ${\mathbb S} \times {\mathbb S}$
solution to the duality problem with given random margins provided  the cost 
function $c$   satisfies
\begin{equation}\label{eq:supertriangle}
   c(x,y)=\sup_{u \in \Lipc} |u(x)-u(y)| \, ,
\end{equation}
where $\Lipc$ is the class of continuous bounded functions $u$ on ${\mathbb S}$
such that $|u(x)-u(y)| \leq c(x,y)$.
As in \cite[Theorem 3.4.1]{cc-prf-valadier04book}, the main tool to prove this result 
is a measurable selection lemma (see Lemma \ref{lem:KRmeas}) 
for an appropriate multifunction. 

Starting from point 1 of Theorem 2.1, we prove in point 2 of Theorem 2.1
that the parametrized Kantorovich--Rubin\v stein theorem
given in \cite[Theorem 3.4.1]{cc-prf-valadier04book} still holds for any 
cost function $c$ satisfying (\ref{eq:supertriangle}). Next, we give in 
Section 3 the application of Theorem 2.1 to the coupling of random variables. In particular,
Corollary \ref{c1} extends Berbee's coupling in the following way:
if $\oap$ is rich enough, and if $c$ is a  mapping satisfying 
(\ref{eq:supertriangle}) such that $\int c(X,x_0) d\pr$ is finite for some
$x_0$ in ${\mathbb S}$, then there exists a random variable $X^*$ distributed as $X$
and independent of ${\mathcal M}$ such that
\begin{equation}\label{eq:couplage}
   \expect{c(X, X^*)}= \Big \| \sup_{f \in \Lipc} \Big | \int f(x) {\pr}_{X|{\cal M}}(dx)- 
\int f(x) {\pr}_{X}(dx) \Big| \, \Big \|_1 \, .
\end{equation}
%be more precise here: let $\oap$ be a probability space,
%let ${\mathcal M}$ be a sub-$\sigma$-algebra of ${\cal A}$ and let  $X$  be a random variable
%with values in ${\mathbb S}$ such that $\int c(X,x_0) d\pr$ is finite for some
%$x_0$ in ${\mathbb S}$. Denote by $\pr_X$ the distribution of $X$ and
%by $\pr_{X|{\mathcal M}}$ a regular conditional distribution of $X$ given ${\mathcal M}$.
If $c(x,y)=\indicatrix_{x \neq y}$ is the discrete metric,  (\ref{eq:couplage}) is exactly Berbee's 
coupling (\ref{ber}). If $c=d$,  
(\ref{eq:couplage}) has been proved in \cite[Corollaire 2.2]{dedecker-prieur04coupl}  
(see also \cite[Section 7.1]{dedecker-prieur03coef}).
For more details on  the coupling property (\ref{eq:couplage}) and its applications, see Section 3.2.

\paragraph{Preliminary notations}

In the sequel, 
For any topological space $\esptopol$, 
we denote by 
$\bor{\esptopol}$ the Borel $\sigma$--algebra of $\esptopol$ and by
%% $\measures{\esptopol}$ the space of finite nonnegative measures 
$\laws{\esptopol}$ the space of probability laws 
on $(\esptopol,\bor{\esptopol})$, 
endowed with the
{narrow topology}, that is, for every mapping $\varphi
:\,{\esptopol}\ra[0,1]$, the mapping 
$\mu\mt\int_\esptopol \varphi\,d\mu$ 
is \lsc if and only if $\varphi$ is \lsc  
%% The subspace of probability laws on $(\esptopol,\bor{\esptopol})$ is
%% denoted by $\laws{\esptopol}$. 

Throughout, $\esp$ is a given completely regular topological space and 
$\oap$ a given probability space. 
Our results are new (at least we hope so) even in the setting of Polish
spaces or simply the real line. 
However they are valid  in much more
general spaces, without significant changes in the proofs. 
{\em The reader who is not interested by this level of 
generality may assume as well in the sequel that all
topological spaces we consider are Polish.} 
On the other hand, we give in appendix some definitions and
references which might be useful for a complete reading.

%%%%%%%%%%%%%%%%%%%%%%%%%%%%%%%%%%%%%%%%%%%%%%%%%%%%%%%%%%
\section{Parametrized Kantorovich--Rubin\v stein theorem}%
\label{section:paramKR}                                  %
%%%%%%%%%%%%%%%%%%%%%%%%%%%%%%%%%%%%%%%%%%%%%%%%%%%%%%%%%%
The results of this section draw inspiration from 
\cite[\S 3.4]{cc-prf-valadier04book}. 
%% We state them in a form which is more general than needed in the
%% sequel, because we feel they are interesting in themselves. 

%% rappels KR 
For any $\mu,\nu\in\laws{\esp}$, let 
$\DistrJointes(\mu,\nu)$ 
be the set of 
probability laws $\pi$ on $(\esp\times\esp,\bor{\esp\times\esp})$ with margins $\mu$ and
$\nu$, that is, $\pi(A\times\esp)=\mu(A)$ and $\pi(\esp\times
A)=\nu(A)$ for every $A\in\bor{\esp}$. 
Let us recall the 

%\begin{theo}\label{theo:KRD}
%\titre{Kantorovich--Rubin\v stein Duality Theorem 
%\cite{levin84topological}, \cite[Theorem 4.6.6]{rachev-ruschendorf98bookI}}
\medskip\noindent
{\bf Kantorovich--Rubin\v stein duality theorem 
\cite{levin84topological}, \cite[Theorem 4.6.6]{rachev-ruschendorf98bookI}}\quad
{\em Assume that $\esp$ is a completely regular pre-Radon space\footnote{
In \cite{levin84topological} and \cite[Theorem 4.6.6]{rachev-ruschendorf98bookI}, the 
%completely regular 
space $\esp$
is assumed to be a universally measurable subset of
some compact space. 
But this amounts to assume that it is completely regular and pre-Radon: 
see \cite[Lemma 4.5.17]{rachev-ruschendorf98bookI} 
and \cite[Corollary 11.8]{gardner-pfeffer84borel-measures}.},  
that is, every finite $\tau$--additive Borel measure on $\esp$ is
inner regular with respect to the compact subsets of $\esp$.  
%see e.g.~\cite[Definition 11.2]{gardner-pfeffer84borel-measures}.  
Let $c :\,\esp\times\esp\ra [0,+\infty]$ be a universally
measurable mapping. For every
$(\mu,\nu)\in\laws{\esp}\times\laws{\esp}$, 
let us denote
\begin{align*}
\KR(\mu,\nu)&:=\inf_{\pi\in\DistrJointes(\mu,\nu)}
               \int_{\esp\times\esp}\cout(\elt,\eltt)\,d\pi(\elt,\eltt),\\
\lipdistc(\mu,\nu)
&:=\sup_{f\in\Lipc}\CCO{\mu(f)-\nu(f)}
\end{align*}
where $
\Lipc=\accol{u\in\Cb{\esp}\tq \foreach{x,y\in\esp}\abs{u(x)-u(y)}\leq\cout(x,y)}$.
Then  the equality 
$
\KR(\mu,\nu)=\lipdistc(\mu,\nu)
$
holds for all $(\mu,\nu)\in\laws{\esp}\times\laws{\esp}$ if and only if (\ref{eq:supertriangle}) holds.
}

\medskip

%\begin{equation*}\label{eq:supertriangle}
%\foreach{x,y\in\esp}
%c(x,y)=\sup_{u\in\Lipc} \abs{u(x)-u(y)}.
%\end{equation*}
Note that, if $c$ satifies (\ref{eq:supertriangle}), it is the supremum of a set of continuous
functions, thus it is l.s.c.
Every continuous metric $\cout$ on $\esp$ satisfies
\eqref{eq:supertriangle} (see 
\cite[Corollary 4.5.7]{rachev-ruschendorf98bookI}), and, if $\esp$ is compact,
every l.s.c.~metric $\cout$  on $\esp$ satisfies
\eqref{eq:supertriangle} (see 
\cite[Remark 4.5.6]{rachev-ruschendorf98bookI}).

%% Definitions KR avec parametre

Now, we denote 
$$\youngs(\esprob,\tribu,\pr;\esp)=\{\mu\in\laws{\esprob\times\esp,\tribu\otimes\bor{\esp}}
\tq \foreach{A\in\tribu}\mu(A\times\esp)=\pr(A)\}.$$
%% by $\youngs$ the set of probability laws 
%% $\mu\in\laws{\esprob\times\esp,\tribu\otimes\bor{\esp}}$ with margin $\pr$ on
%% $\esprob$, that is,  $\mu(A\times\esp)=\pr(A)$ for every
%% $A\in\tribu$.
When no confusion can arise, we omit some part of the information, and
use notations such as $\youngs(\tribu)$ or simply $\youngs$ 
(same remark for the set $\youngsint{\cout}(\esprob,\tribu,\pr;\esp)$
defined below). 
If $\esp$ is a Radon space, every $\mu\in\youngs$
is {\em disintegrable}, that is, there exists a (unique, up to
$\pr$-a.e.~equality) ${\mathcal A}^*_\mu$-measurable mapping
$\elprob\mt\mu_\elprob$, $\esprob\ra\laws{\esp}$, such that 
$$\mu(f)=\int_\esprob\int_\esp
f(\elprob,\elt)\,d\mu_\elprob(\elt)\,d\pr(\elprob)$$ 
for every measurable $f :\,\esprob\times\esp\ra[0,+\infty]$ 
(see \cite{valadier73desi}). If furthermore the compact subsets of $\esp$ are metrizable, the
mapping $\elprob\mapsto\mu_\elprob$ can be chosen $\tribu$-measurable, 
see the Appendix.

Let $\cout$ satisfy  \eqref{eq:supertriangle}. 
We denote $$\youngsint{\cout}(\esprob,\tribu,\pr;\esp)=\{\mu\in\youngs\tq 
\int_{\esprob\times\esp} \cout(\elt,\elt_0)\,d\mu(\elprob,\elt)<+\infty\}$$
where $\elt_0$ is some fixed element of $\esp$ (this definition is
independent of the choice of $\elt_0$). 
  For any $\mu,\nu\in\youngs$, let 
  $\paramDistrJointes(\mu,\nu)$ be the set of 
  probability laws $\pi$ on $\esprob\times\esp\times\esp$ such that
  $\pi(.\times . \times\esp)=\mu$ and 
  $\pi(.\times \esp \times .)=\nu$. 
We now define the parametrized versions of $\KR$ and $\lipdistc$. 
  Set, for $\mu,\nu\in\youngsint{\cout}$, 
  $$
  \paramKR(\mu,\nu)
  =\inf_{\pi\in\paramDistrJointes(\mu,\nu)}
  \int_{\esprob\times\esp\times\esp} 
               \cout(\elt,\eltt)\,d\pi(\elprob,\elt,\eltt).
  $$
Let also $\paramLipc$ denote the set of measurable integrands 
$f :\,\esprob\times\esp\ra \R$ such that $f(\elprob,.)\in\Lipc$
for every $\elprob\in\esprob$. 
We denote 
$$\paramlipdistc(\mu,\nu)
=\sup_{f\in\paramLipc}\CCO{\mu(f)-\nu(f)}.$$

%%%%%%%%%%%%%% PRAMETRIZED KANTOROVICH - RUBINSTEIN THEOREM 
\begin{theo}\label{theo:KW}
\titre{Parametrized Kantorovich--Rubin\v stein theorem}
Assume that $\esp$ is a completely regular Radon space and that the
compact subsets of $\esp$ are metrizable (e.g.~$\esp$ is a regular
Suslin space). 
Let $\cout :\,\esp\times\esp\ra[0,+\infty[$ 
satisfy \eqref{eq:supertriangle}. 
Let $\mu,\nu\in\youngsint{\cout}$ and let $\elprob\mt\mu_\elprob$ and
$\elprob\mt\nu_\elprob$ be disintegrations of $\mu$ and $\nu$
respectively.  
\begin{enumerate}
\item Let $G :\, \elprob\mt
  \KR(\mu_\elprob,\nu_\elprob)=\lipdistc(\mu_\elprob,\nu_\elprob)$ and
  let  $\complet{\tribu}$ be the universal completion of
  $\tribu$. There exists an $\complet{\tribu}$--measurable mapping
  $\elprob\mt\lambda_\elprob$ from $\esprob$ to ${\cal P}(\esp\times \esp)$ such
  that $\lambda_\elprob$ belongs to
  $\DistrJointes(\mu_\elprob,\nu_\elprob)$ and 
\begin{equation*}
G(\elprob)=\int_{\esp\times\esp}\cout(\elt,\eltt)\,d\lambda_\elprob(\elt,\eltt).
\end{equation*}
\item The following equalities hold:
\begin{equation*}
\paramKR(\mu,\nu)
=\int_{\esprob\times\esp\times\esp}\cout(\elt,\eltt)\,d\lambda(\elprob,\elt,\eltt)
=\paramlipdistc(\mu,\nu),
\end{equation*}
where $\lambda$ is the element of
$\youngs(\esprob,\tribu,\pr;\esp\times\esp)$ defined by 
$\lambda(A\times B\times C)=\int_A \lambda_\elprob(B\times
C)\,d\pr(\elprob)$ for any $A$ in $\tribu$, $B$ and $C$ in
$\bor{\esp}$. In particular, $\lambda$ belongs to
$\paramDistrJointes(\mu,\nu)$, and the infimum in the definition of
$\paramKR(\mu,\nu)$ is attained for this $\lambda$. 
\end{enumerate}
\end{theo}
%%%%%%%%%%%%%
Let us first prove the following lemma. 
The set of compact subsets of a topological space $\esptopol$ is denoted by
$\compacta(\esptopol)$. 

%%%%%%%%%%%%%%%%%% Lem selection mesurable K-R %%
\begin{lem}\label{lem:KRmeas}
\titre{A measurable selection lemma}
Assume that $\esp$ is a Suslin space. 
Let $\cout :\,\esp\times\esp\ra[0,+\infty]$ be an l.s.c.~mapping. 
Let $\bor{}^*$ be the universal completion of the $\sigma$--algebra 
$\bor{\laws{\esp}\times\laws{\esp}}$.  
For any $\mu,\nu\in\laws{\esp}$, let 
$$\niveaubas(\mu,\nu)=\inf_{\pi\in\DistrJointes(\mu,\nu)}
\int \cout(\elt,\eltt)\,d\pi(\elt,\eltt)\in[0,+\infty].$$
The function $\niveaubas$ is $\bor{}^*$--measurable. Furthermore, 
the multifunction  
$$
K :\,
\left\{\begin{array}{lcl}
\laws{\esp}\times\laws{\esp}&\ra&\compacta\CCO{\laws{\esp\times\esp}}\\
(\mu,\nu)&\mt&\accol{\pi\in \DistrJointes(\mu,\nu) 
\tq \int\cout(\elt,\eltt)\,d\pi(\elt,\eltt)
=\niveaubas(\mu,\nu)}
\end{array}\right.
$$
has a $\bor{}^*$--measurable selection,
that is, there exists a $\bor{}^*$--measurable mapping 
$\lambda :\, (\mu,\nu)\mapsto \lambda_{\mu,\nu}$ defined on
$\laws{\esp}\times\laws{\esp}$ with values in
$\compacta\CCO{\laws{\esp\times\esp}}$, such that 
$\lambda_{\mu,\nu}\in K(\mu,\nu)$ for all $\mu,\nu\in\laws{\esp}$.
%%$\elprob\mt\pi_\elprob$. 
\end{lem}
%%%%%%%%%
\preuve 
Observe first that 
the mapping $\niveaubas$ can be defined as 
$$\niveaubas :\,
(\mu,\nu)\mt \inf\accol{\psi(\pi)\tq \pi\in \DistrJointes(\mu,\nu)},$$ 
with 
$$\psi :\, \left\{\begin{array}{lcl}
\laws{\esp\times\esp}&\ra&[0,+\infty]\\
\pi&\mt&\int_{\esp\times\esp}\cout(\elt,\eltt)\,d\pi(\elt,\eltt). 
\end{array}\right.$$
The mapping $\psi$ is \lsc because it is
the supremum of the \lsc mappings $\pi\mt\pi(\cout\wedge n)$,
$n\in\N$ 
(if $\cout$ is bounded and continuous, $\psi$ is continuous). 
Furthermore, we have 
$\DistrJointes=\Phi^{-1}$, where $\Phi$ is the continuous mapping 
$$
\Phi :\,
\left\{\begin{array}{lcl}
\laws{\esp\times\esp}&\ra&\laws{\esp}\times\laws{\esp}\\
\lambda&\mt&(\lambda(.\times \esp),\lambda(\esp\times .))
\end{array}\right.
$$
(recall that $D(\mu,\nu)$ is the set of 
probability laws $\pi$ on $\esp\times\esp$ with margins $\mu$ and
$\nu$). 
Therefore, the graph $\graphe{\DistrJointes}$ 
of $\DistrJointes$ is a closed subset of the Suslin space 
$\X=\bigl( \laws{\esp}\allowbreak\times\laws{\esp} \bigr)\times\laws{\esp\times\esp}$.  
%%see page \pageref{alsosuslin}
%%thus it is measurable. 
Thus, for every $\alpha\in\R$, the set 
 $$\accol{((\mu,\nu),\pi)\in\graphe{\DistrJointes}\tq
\psi(\pi)<\alpha}$$
is a Suslin subset of $\X$. We thus have, by  the \projtheo\ (see
\cite[Lemma III.39]{castaing-valadier77book}),  
\begin{multline*}
\foreach{\alpha\in\R}\accol{(\mu,\nu)\tq \niveaubas(\mu,\nu)<\alpha}\\
=\proj{\laws{\esp}\times\laws{\esp}}\accol{((\mu,\nu),\pi)\in\graphe{\DistrJointes}\tq
\psi(\pi)<\alpha}
\in\bor{}^*.
\end{multline*}

%%%%%%%%%%%%%%%%%%%%%%%%%%%%%%%%%

Now, for each $(\mu,\nu)\in\laws{\esp}\times\laws{\esp}$, we have 
$$K(\mu,\nu)=
%%\defmap{\laws{\esp}\times\laws{\esp}}{\compacta\CCO{\laws{\esp\times\esp}}}%
{\accol{\pi\in \DistrJointes(\mu,\nu)\tq \psi(\pi)=\niveaubas(\mu,\nu)}.}$$
The multifunction $K$ has nonempty compact values because $\DistrJointes$ has nonempty
compact values and $\psi$ is l.s.c. 
Let 
$$F
 :\,\left\{\begin{array}{lcl}
\CCO{\laws{\esp}\times\laws{\esp}}\times\laws{\esp\times\esp}&\ra&\R\\
{((\mu,\nu),\pi)}&\mt&{\psi(\pi)-\niveaubas(\mu,\nu).}
\end{array}\right.$$
The mapping $F$ is $\bor{}^*\otimes\bor{\laws{\esp\times\esp}}$--measurable. 
Furthermore, the graph of $K$ is 
\begin{align*}
\graphe{K}
&=\accol{((\mu,\nu),\pi)\tq \mu=\pi(.\times\esp),\ \nu=\pi(\esp\times.),\ 
                            F((\mu,\nu),\pi)=0}\\
&=\graphe{\DistrJointes}\cap F^{-1}(0)\\
&\in\bor{}^*\otimes\bor{\laws{\esp\times\esp}}.
\end{align*}
As $\esp$ is Suslin, this proves that $K$ is 
$\bor{}^*$--measurable (see \cite[Theorem III.22]{castaing-valadier77book}).  %(see Remark \ref{suslin}). 
Thus $K$ has a $\bor{}^*$--measurable selection. 
\fin
%%%%%

%% \comm{Dans le lemme 3.4.2 de \cite{cc-prf-valadier04book}, on se 
%% ramenait au cas polonais car $\esp$ n'\'etait pas n\'ecessairement
%% souslinien. 
%% En fait, on doit pouvoir prendre $\esp$ cosmique et se ramener au cas
%% souslinien... 
%% l'int\'er\^et des espaces cosmiques est qu'ils englobent \`a la fois
%% les espaces m\'etriques s\'eparables et les espaces sousliniens. }

%%%%%%%%%%%%%%%%%%%%%%%%%%%%%%%%%
\preuvof{Theorem \ref{theo:KW}} %
By the Radon property, the probability measures $\mu(\esprob\times.)$
and $\nu(\esprob\times.)$ are tight, that is, for every integer $n\geq 1$, 
there exists a compact subset $K_n$ of $\esp$ such that 
$\mu(\esprob\times (\esp\setminus K_n))\leq 1/n$ and 
$\nu(\esprob\times (\esp\setminus K_n))\leq 1/n$. 
Now, we can clearly replace $\esp$ in the statements of Theorem
\ref{theo:KW} by the smaller space $\cup_{n\geq 1}K_n$. 
But $\cup_{n\geq 1}K_n$ is Suslin (and even Lusin), so we can assume without loss of
generality that $\esp$ is a regular Suslin space.

%% Furthermore, from the usual 
%% Kantorovich--Rubin\v stein Theorem \cite[Theorem
%% 4.6.6]{rachev-ruschendorf98bookI}, 
%% we have, for every $\elprob\in\esprob$, 
%% \begin{equation}\label{eq:coherent}
%% \inf_{\pi\in\DistrJointes(\mu_\elprob,\nu_\elprob)}
%%                {\int_{\esp\times\esp}
%%                  \cout(\elt,\eltt)}\,d\pi(\elt,\eltt)
%% =\sup_{g\in\Lipc}\CCO{\mu_\elprob(g)-\nu_\elprob(g)}, 
%% \end{equation}
%% which proves that the definition of $G$ is coherent and proves in the
%% same time 
%% the second equality in \eqref{eq:paramKR}. 
We easily have  
\begin{align}
\paramlipdistc(\mu,\nu)
&=\sup_{f\in\paramLipc}\int_\esprob
 \int_\esp\int_\esp \CCO{f(\elprob,\elt)-f(\elprob,\eltt)}\,d\mu_\elprob(\elt)\,d\nu_\elprob(\eltt)
 \,d\pr(\elprob) \notag\\
&\leq \int_\esprob
 \int_\esp\int_\esp \cout(\elt,\eltt)\,d\mu_\elprob(\elt)\,d\nu_\elprob(\eltt)
 \,d\pr(\elprob) \notag\\
&\leq \paramKR(\mu,\nu). \label{eq:inegalite-facile}
\end{align}
So, to prove Theorem \ref{theo:KW}, we only need to prove that
$\paramKR(\mu,\nu)\leq \paramlipdistc(\mu,\nu)$ and that the minimum
in the definition of $\paramKR(\mu,\nu)$ is attained.

%{\em First step.}  
%By \eqref{eq:supertriangle}, the function $\cout$ id \lsc
Using the notations of Lemma \ref{lem:KRmeas}, we have 
$G(\elprob)=\niveaubas(\mu_\elprob,\nu_\elprob),$ 
thus $G$ is $\complet{\tribu}$--measurable   
(indeed, the mapping $\elprob\mt(\mu_\elprob,\nu_\elprob)$ is
               measurable for $\complet{\tribu}$ and $\bor{}^*$
               because it is measurable for $\tribu$ and  
               $\bor{\laws{\esp}\times\laws{\esp}}$). 
From Lemma \ref{lem:KRmeas}, 
the multifunction  $\elprob\mt
\DistrJointes(\mu_\elprob,\nu_\elprob)$ has an $\complet{\tribu}$--measurable selection 
$\elprob\mt\lambda_\elprob$ such that, for every $\elprob\in\esprob,$
 $G(\elprob)=\int_{\esp\times\esp}
                 \cout(\elt,\eltt)\,d\lambda_\elprob(\elt,\eltt).$ 
We thus have 
\begin{equation}\label{eq:voila-le-min}
\paramKR(\mu,\nu)\leq 
\int_{\esprob\times\esp\times\esp} 
                \cout(\elt,\eltt)\,d\lambda(\elprob,\elt,\eltt)
%%\allowbreak
=\int_\esprob G(\elprob)\,d\pr(\elprob).
\end{equation}
Furthermore, since $\mu,\nu\in\youngsint{\cout}$, we have 
$G(\elprob)<+\infty$ a.e. 
Let $\esprob_0$ be the almost sure set on which
$G(\elprob)<+\infty$. Fix an element $x_0$ in $\esp$. 
%% From the usual 
%% Kantorovich--Rubin\v stein Theorem \cite[Theorem
%% 4.6.6]{rachev-ruschendorf98bookI}, 
We have, for every $\elprob\in\esprob_0$, 
\begin{equation*}%\label{eq:gphK}
G(\elprob)=
\sup_{g\in\Lipc}\CCO{\mu_\elprob(g)-\nu_\elprob(g)}
=\sup_{\scriptstyle g\in\Lipc,\ 
                        \scriptstyle g(\elt_0)=0}
                              \CCO{\mu_\elprob(g)-\nu_\elprob(g)}.
\end{equation*}
Let $\epsilon>0$. 
Let $\widetilde{\mu}$ and $\widetilde{\nu}$ be the finite measures on
$\esp$ 
defined by 
$$\widetilde{\mu}(B)
=\int_{\esprob\times  B}\cout(\elt_0,\elt)\,d\mu(\elprob,\elt)
\quad \text{ and }\quad 
\widetilde{\nu}(B)
=\int_{\esprob\times  B}\cout(\elt_0,\elt)\,d\nu(\elprob,\elt)
$$ 
for any $B\in\bor{\esp}$. 
Let $\espp$ be a compact subset of $\esp$ containing $\elt_0$ such that 
$\widetilde{\mu}({\esp\setminus\espp})\leq \epsilon$ and 
$\widetilde{\nu}({\esp\setminus\espp})\leq \epsilon$.  
For any $f\in\paramLipc$, we have 
\begin{multline}
\abs{\int_\esprob (\mu_\elprob-\nu_\elprob)(f(\elprob,.))\,d\pr(\elprob)
-\int_\esprob
  (\mu_\elprob-\nu_\elprob)(f(\elprob,.)\un{\espp})\,d\pr(\elprob)
}\\
= \abs{\int_\esprob
  (\mu_\elprob-\nu_\elprob)(f(\elprob,.)\un{{\esp\setminus\espp}})\,d\pr(\elprob)}%%\notag \\ &
\leq 2\epsilon. \label{gmoinsgtz}
\end{multline}
Set, for all $\elprob\in\esprob_0$,  
$$G'(\elprob)=\sup_{\begin{array}{c}\scriptstyle g\in\Lipc, \ 
                       \scriptstyle g(\elt_0)=0\end{array}}
       (\mu_\elprob-\nu_\elprob)(g\un{\espp}).$$ 
We thus have 
\begin{equation} \label{GmoinsGprime}
\abs{\int_{\esprob_0} G\,d\pr - \int_{\esprob_0} G'\,d\pr}\leq 2\epsilon. 
\end{equation}
Let $\Lipc\restr{\espp}$ denote the set of restrictions to $\espp$ of
elements of $\Lipc$. 
The set $\espp$ is metrizable, thus $\Cb{\espp}$ (endowed with the
topology of uniform convergence) is metrizable separable, thus its subspace 
$\Lipc\restr{\espp}$ is also metrizable separable. We can thus find a dense
countable subset $\dense=\{u_n\tq n\in\N\}$ of $\Lipc$ 
for the seminorm
$\norm{u}_{\Cb{\espp}}:=\sup_{\elt\in\espp}\abs{u(\elt)}$. 
Set, for all $(\elprob,\elt)\in\esprob_0\times \esp$,   
\begin{align*}
&N(\elprob)=\min\bigl\{  n\in\N\tq \int_\esp
u_n(\elt)\,d(\mu_\elprob-\nu_\elprob)(\elt)\geq
%\lipdistc(\mu_\elprob,\nu_\elprob)-\epsilon
G'(\elprob)-\epsilon
                  \bigr\},\\ 
&f(\elprob,\elt)=u_{N(\elprob)}(\elt).\tag*{and}
\end{align*}
We then have, using \eqref{gmoinsgtz} and \eqref{GmoinsGprime},    
\begin{align*}
\paramlipdistc(\mu,\nu)
\geq \int_{\esprob_0\times \esp} f\,d(\mu-\nu)
&\geq \int_{\esprob_0\times \espp} f\,d(\mu-\nu)-2\epsilon \\
&\geq \int_{\esprob_0} G'\,d\pr - 3\epsilon 
\geq \int_{\esprob_0} G\,d\pr - 5\epsilon.
\end{align*}
Thus, in view of 
\eqref{eq:inegalite-facile} and\eqref{eq:voila-le-min}, 
$$\paramKR(\mu,\nu)=\int_{\esprob\times\esp\times\esp} 
                \cout(\elt,\eltt)\,d\lambda(\elprob,\elt,\eltt)
=\paramlipdistc(\mu,\nu).$$
%which proves Theorem \ref{theo:KW}.
\fin
%%%%%

%%%%%%%%%%%%%%%%%%%%%%%%%%%%%%%%%%%%%%%%%%%%%%%%%%%%%%%%%
\section{Application: coupling for the minimal distance}
%%%%%%%%%%%%%%%%%%%%%%%%%%%%%%%%%%%%%%%%%%%%%%%%%%%%%%%%%

%\subsection{Main coupling result}

In this section   ${\mathbb S}$  is  a completely regular Radon space with metrizable
compact subsets,  
$c:{\mathbb  S}\times {\mathbb S}\rightarrow [0, +\infty]$ is 
a  mapping satisfying \eqref{eq:supertriangle} and
${\cal M}$ is a sub-$\sigma$-algebra of ${\cal A}$. Let 
$X$ be a random variable with values in   $\mathbb{S}$, let $\pr_X$ be the distribution of 
$X$, and let    $\pr_{X|{\cal M}}$ be a regular
conditional distribution of $X$ given ${\cal M}$ (see Section  \ref{rap} for the  existence).
We assume that $\int c(x, x_0) \pr_X(dx)$ is finite for some (and therefore any) $x_0$ in
${\mathbb S}$ (which means exactly that the unique measure of $\youngs ({\cal M})$
with disintegration $\pr_{X|{\cal M}}(\cdot, \omega)$ belongs to 
$\youngsint{\cout}({\cal M})$).

\begin{theo}[general coupling theorem]\label{p2} Assume that $\Omega$ is rich enough, that is, there exists
a random variable  $U$ from $(\Omega, \mathcal{A})$ to $([0,1], \mathcal{B}([0,1]))$,
independent of  $\sigma(X) \vee {\cal M}$ and uniformly distributed
over $[0,1]$. Let $Q$ be any element of $\youngsint{\cout} ({\cal M})$. 
%such that
%$\omega \mapsto Q_\omega(A)$ is ${\cal M}$-measurable for any $A$ of ${\cal B}_{\mathbb S}$.
%and $P$  be the measure of ${\youngs}({\cal M})$ whose disintegrated measure is
%$P_\omega( \cdot)={\pr}_{X |{\cal M}}(\cdot, \omega)$.
%If  $P$ and $Q$ belong to  $\youngsint{\cout}({\cal M})$, 
%Soit $\mu$ la probabilit\'e r\'eguli\`ere de la Proposition \ref{p1}
There exists a  $\sigma(U) \vee \sigma(X)\vee {\cal M}$-measurable random variable $Y$, such that
$Q_\point$ is a regular conditional probability of $Y$ given ${\cal M}$, and 
\begin{equation}\label{(3)}
   \condexpect{c(X,Y)}{\cal M}=
\sup_{f \in \Lipc} \Big | \int f(x) {\mathbb P}_{X|{\cal M}}(dx)- 
\int f(x) Q_\point(dx) \Big| \quad \text{${\pr}$-a.s.} \, .
\end{equation}
\end{theo}
%On dispose pour cela d'une variable $U$ de ${\Omega}$ dans $[0,1]$, dont la loi
%est uniforme sur $[0,1]$, et ind\'ependante de $ \sigma(Z,X)$ (c'est
%\`a dire qu'on suppose que $\Omega$ est ``assez riche'').

\medskip

\preuve 
%Since ${\mathbb S}$ is Lusin, it is Suslin.
%As in the proof of Theorem \ref{theo:KW}, we  assume without loss that ${\mathbb S}$ is Lusin regular.
We apply  Theorem \ref{theo:KW} to the probability space $(\Omega, {\cal M}, \pr)$
and to the disintegrated measures $\mu_\omega(\cdot)=\pr_{X|{\cal M}}(\cdot, \omega)$ and $\nu_\omega=Q_\omega$.
As in the proof of Theorem \ref{theo:KW}, we  assume without loss of generality that ${\mathbb S}$ is Lusin regular.
 From  point 1 of Theorem \ref{theo:KW} we infer that there exists a
mapping $\omega \mapsto \lambda_\omega$
from $\Omega$ to ${\cal P}({\mathbb S}\times {\mathbb S})$, measurable for 
${\cal M}^*$ and ${\cal B}_{{\cal P}({\mathbb S}\times {\mathbb S})}$, such that  $\lambda_\omega$ belongs to 
$D({\pr}_{X|{\cal M}}(\cdot, \omega), Q_\omega)$
and $G(\omega)=\int_{\mathbb S \times \mathbb S} c(x,y) \lambda_\omega(dx,dy)$.

On the measurable space $({\mathbb M}, {\mathcal T})=
(\Omega \times {\mathbb S} \times {\mathbb S},{\mathcal M}^* 
\otimes {\mathcal B}_{\mathbb S} \otimes {\mathcal B}_{\mathbb S})$ 
we put the probability
$$
     \pi(A\times B \times C) = \int_A \lambda_\omega(B \times C) \pr(d\omega)\, .
$$
If $I=(I_1, I_2, I_3)$ is the identity  on  ${\mathbb M}$, we see that 
%$\omega \mapsto \lambda_\omega$ is 
a regular conditional distribution of $(I_2, I_3)$ given $I_1$
is given by $\pr_{(I_2, I_3)|I_1=\omega}=\lambda_\omega$. 
Since  ${\pr}_{X|{\cal M}}(\cdot, \omega)$ is the first margin of $\lambda_\omega$, 
 a regular conditional probability of $I_2$ given $I_1$ is given by 
$\pr_{I_2|I_1=\omega}(\cdot)={\pr}_{X|{\cal M}}(\cdot, \omega)$. Let
$ \lambda_{\omega,x}=\pr_{I_3|I_1=\omega, I_2=x}$ be a regular conditional
distribution of $I_3$ given $(I_1, I_2)$,  so that 
$(\omega, x) \mapsto \lambda_{\omega, x}$  is  measurable for
${\cal M}^*\otimes {\cal B}_{\mathbb  S}$ and 
${\cal B}_{{\cal P}({\mathbb S})}$.  
From the unicity (up to $\pr$-a.s. equality) of regular conditional probabilities,  it follows that 
%Since  ${\pr}_{X|{\cal M}}(\cdot, \omega)$ is the first margin of $\lambda_\omega$, we infer from the 
%disintegration theorem (see \cite{valadier73desi}) that there exists
%a mapping $(\omega, x) \mapsto \lambda_{\omega, x}$ from
%$\Omega \times {\mathbb S}$ to ${\cal P}({\mathbb S})$, measurable for
%${\cal M}^*\otimes {\cal B}_{\mathbb  S}$ and 
%${\cal B}_{{\cal P}({\mathbb S})}$, such that
\begin{equation}
   \lambda_\omega(B \times C)= \int_B \lambda_{\omega,x}(C){\pr}_{X|{\cal M}}(dx, \omega) 
 \quad \text{$\pr$-a.s.} \, .
\end{equation}
%Rappelons que sur 
% $(S, \mathcal{B}(S))=(\Omega \times \mathcal{X} \times \mathcal{X}, 
%{\cal M}^* \otimes \mathcal{B}(\mathcal{X})\otimes \mathcal{B}(\mathcal{X}))$
%on d\'efinit
%\begin{equation}\label{(4)}
%    \lambda(A\times B\times C)= \int_A \lambda_\omega(B\times C) {\mathbb P}(d\omega) \, .
%\end{equation}
%On note $I=(I_1, I_2, I_3)$ l'identit\'e de $S$ dans $S$.
% de sorte
%que la loi de $I$ est $Q$. 

Assume that we can find
a random variable 
$\tilde{Y}$ from $\Omega$ to $\mathbb S$, measurable for $\sigma(U) \vee \sigma(X)\vee {\cal M}^*$
and $\mathcal{B}_{\mathbb{S}}$, such that
${\pr}_{\tilde{Y}|\sigma(X)\vee {\cal M}^*}(\cdot, \omega)= \lambda_{\omega, X(\omega)}(\cdot)$.
Since $\omega \mapsto {\pr}_{X |{\cal M}}(\cdot, \omega)$ is measurable for ${\cal M}^*$ and
${\cal B}_{{\cal P}({\mathbb S})}$, one can check that
${\pr}_{X |{\cal M}}$ is a regular conditional probability of
$X$ given ${\cal M}^*$. For $A$ in ${\cal M}^*$, $B$ and $C$ in
${\cal B}_{\mathbb S}$, we thus have
\begin{eqnarray*}
\expect{\indicatrix_A \indicatrix_{X \in B} \indicatrix_{\tilde{Y} \in C}}
&=&
\expect{\indicatrix_A\condexpect{\indicatrix_{X \in B} 
\condexpect{\indicatrix_{\tilde{Y} \in C}}{\sigma(X)\vee{\cal M}^*}}{{\cal M}^*}}\\
&=&
\int_A \Big(\int_B \lambda_{\omega,x}(C) {\pr}_{X|{\cal M}}(dx, \omega)\Big) {\pr}(d\omega) \\
&=& 
\int_A \lambda_\omega(B \times C) {\pr}(d\omega) \, .
\end{eqnarray*}
We infer that $\lambda_\omega$ is a regular conditional probability of 
$(X, \tilde{Y})$ given ${\cal M}^*$.
% ce qui entra\^{\i}ne (\ref{(3)}). 
By definition of  $\lambda_\omega$, we obtain that
\begin{equation}\label{etoile}
   \condexpect{c(X,\tilde{Y})}{{\cal M}^*}= 
\sup_{f \in \Lipc} \Big | \int f(x) {\pr}_{X|{\cal M}}(dx)- 
\int f(x) Q_\point(dx) \Big| \quad \text{ ${\pr}$-a.s.} \, .
\end{equation}
%{\mathbb P}_{I_3|I_2=X(\omega), I_1=\omega}(\cdot)$.
%La distribution de $(I_1, X, \tilde{Y})$ est alors donn\'ee par: 
%pour tout $(A, B, C)$ dans 
%$\mathcal{B}(\mathcal{X}) \times \mathcal{B}(\mathcal{X}) \times \mathcal{B}(\mathcal{Z})$,
%$$
%    {\mathbb P}_{I_1, X, \tilde{Y}}(A \times B \times C)=
%\int_{A \times C} {\mathbb P}_{\tilde{Y}|X=x, Z=z}(B) {\mathbb P}_{X, Z} (dx, dz) 
%=
%\int_{A} \Big (\int_B{\mathbb P}_{I_3|I_2=x, I_1=\omega}(C) 
% {\mathbb P}_{X|{\cal M}^*}(dx, \omega) \Big ) {\mathbb P}(d\omega)\, .
%$$
%Comme ${\mathbb P}_{X | {\cal M}^*}( \cdot, \omega)=\mu_\omega(\cdot)$, on en 
%d\'eduit que 
%Gr\^ace au point {\it (ii)} de la proposition \ref{p1}, on a  ${\mathbb P}_{X, Z} 
%est donn\'ee par: pour tout $(A, B)$ dans
%$(\mathcal{B}(\mathcal{Z})\times\mathcal{B}(\mathcal{X}))$,
%$$
%   {\mathbb P}_{Z, X}(A\times B)=\int_A {\mathbb P}_{X|Z=z}(B) {\mathbb P}_{Z}(dz) \, .
%$$
%En utilisant 2, on obtient
%$$
%    {\mathbb P}_{Z, X}(A\times B)=\int_A \mu( B \times \mathcal{X},z) {\mathbb P}_Z(dz)=
%    Q(A \times B \times \mathcal{X})={\mathbb P}_{I_1, I_2}(A \times B) \, .
%$$
%={\mathbb P}_{I_1, I_3}$, et donc
%$ {\mathbb P}_{X, \tilde{Y}, Z}=\nu$. 
%On d\'eduit  de (\ref{(4)}) que $\lambda_\omega$ est une distribution conditionnelle de $(X, \tilde{Y})$ sachant
%${\cal M}^*$
Since ${\mathbb S}$ is Lusin, it  is standard Borel (see Section \ref{rap}).
Applying Lemma \ref{lem:modification}, 
there exists a $\sigma(U) \vee \sigma(X)\vee {\cal M}$-measurable modification
$Y$ of  $\tilde{Y}$, so that
(\ref{etoile})
still holds for  ${\mathbb E}(c(X,Y)|{\cal M}^*)$. We obtain (\ref{(3)}) by noting that 
$\condexpect{c(X,Y)}{{\cal M}^*}=\condexpect{c(X,Y)}{{\cal M}}$ ${\pr}$-a.s.
%Pour obtenir (\ref{(3)}), on prend une modification 
%il suffit de prendre une modification $\sigma(U) \vee \sigma(X)\vee {\cal M}$-mesurable.1
%$\tilde{Y}$ de $\tilde{Y}$.
%$$
%  {\mathbb P}_{X,\tilde{Y}}(A)=Q(\Omega \times A)=\int \mu(A, \omega) {\mathbb P}(d\omega)
%$$
%et on obtient la loi jointe donn\'ee dans (2). Il reste \`a v\'erifier que 
%${\mathbb P}_{\tilde{Y}|I_1}={\mathbb P}_X$. Mais pour $(A, B)$ dans $(\mathcal{M}\times\mathcal{B}(\mathcal{X}))$, on a, 
%en vertu de 3
%$$
%   {\mathbb E}({\mathbb I}_{I_1 \in A} {\mathbb I}_{\tilde{Y} \in B})=\int_A \mu(\mathcal{X}\times B) {\mathbb P}(d\omega) 
%    = {\mathbb P}(A) {\mathbb P}_X(B) \, ,
%$$ 
%ce qui permet de conlure.

It remains to build $\tilde{Y}$. Since ${\mathbb S}$  is standard Borel, there exists a one to one map
$f$ from $\mathbb S$ to a Borel subset of  $[0,1]$, such that $f$ and $f^{-1}$ are measurable
for ${\cal B}([0,1])$ and ${\cal B}_{\mathbb S}$. 
%Par cons\'equent $f(I_2)$ est une variable al\'eatoire  et  
%${\mathbb P}_{f(I_2)|I_1, I_3}(A)={\mathbb P}_{I_2|I_1, I_3}(f^{-1}(A))$.
Define $F(t,  \omega)=\lambda_{\omega, X(\omega)}(f^{-1}(]-\infty, t]))$. The map
$F(\cdot,  \omega)$ is a  distribution function with  c\`adl\`ag inverse
$F^{-1}(\cdot,  \omega)$. One can see that the map
$(u,\omega) \rightarrow F^{-1}(u,\omega)$ is
$\mathcal{B}([0,1])\otimes {\cal M}^*\vee \sigma(X)$-measurable. Let 
$T(\omega)=F^{-1}(U(\omega), \omega)$ and 
$\tilde{Y}=f^{-1}(T)$. It remains to see that
${\pr}_{\tilde{Y}|\sigma(X) \vee {\cal M}^*}(\cdot, \omega)= \lambda_{\omega, X(\omega)}(\cdot)$.
For any $A$ in ${\cal M}^*$, $B$  in 
${\cal B}_\mathbb S$ and $t$ in ${\mathbb R}$, we have 
$$
   \expect{ \indicatrix_{ A}\indicatrix_{X\in B}\indicatrix_{\tilde{Y} \in f^{-1}(]-\infty, t])}}= 
   \int_A \indicatrix_{X(\omega) \in B} 
   \indicatrix_{U(\omega) \leq F(t, \omega)} {\pr}(d\omega).
$$
Since $U$ is independent of $\sigma(X) \vee {\cal M}$, it is also 
independent of $\sigma(X) \vee {\cal M}^*$. Hence
\begin{eqnarray*}
\expect{\indicatrix_{A} \indicatrix_{X \in B} \indicatrix_{\tilde{Y} \in f^{-1}(]-\infty, t])}} &=&
\int_A \indicatrix_{X(\omega) \in B} 
F(t, \omega)\pr (d\omega)\\
%&=&
%\int {\mathbb I}_{X(\omega) \in A} {\mathbb I}_{Z(\omega) \in B}
%{\mathbb P}_{I_2|I_1=X(\omega), I_3=Z(\omega)}(f^{-1}(]-\inft\tilde{\tilde{Y}}, t])){\mathbb P}(d\omega) \\
&=&
\int_A \indicatrix_{X(\omega) \in B} 
\lambda_{\omega, X(\omega)}(f^{-1}(]-\infty, t])) \pr (d\omega).
\end{eqnarray*}
Since $\{f^{-1}(]-\infty, t]) , t \in [0,1] \}$ is a separating class,  the result follows.
\fin

\paragraph{Coupling and dependence coefficients}
Define the coefficient
\begin{equation}\label{(5)}
   \tau_c(\mathcal{M}, X)=\Big \| \sup_{f \in \Lipc} \Big | \int f(x) {\pr}_{X|{\cal M}}(dx)- 
\int f(x) {\pr}_{X}(dx) \Big| \, \Big \|_1 \, .
\end{equation}
%La version uniforme de ce  coefficient \`a \'et\'e introduite par Rio dans \cite{Ri}.
If $\Lipc$ is a separating class, this coefficient measures the dependence between $\mathcal{M}$ and
$X$ ($\tau_c(\mathcal{M}, X)=0$ if and only if $X$ is independent of ${\cal M}$). From point 2
of Theorem  \ref{theo:KW}, we see that an equivalent definition is
$$
 \tau_c(\mathcal{M}, X)=\sup_{f \in \Lipcm}
 \int f(\omega, X(\omega)) {\pr}(d\omega) -
 \int \Big (\int f(\omega, x) {\pr}_X(dx) \Big ) {\pr}(d\omega) \, .
$$
where $\Lipcm$ is the set of integrands  $f$
from $\Omega \times {\mathbb S} \rightarrow {\mathbb R}$, measurable for ${\cal M}\otimes {\cal B}_{\mathbb S}$, 
such that
$f(\omega,.)$ belongs to $ \Lipc$ for any $\omega \in \Omega$.

Let $c(x,y)= \indicatrix_{x \neq y}$ be the discrete metric and let $\|\cdot\|_v$ be the variation norm. 
From the Riesz-Alexandroff representation theorem (see \cite[Theorem 5.1]{wheeler83survey}), we infer that 
for any $(\mu, \nu)$ in ${\cal P}(\mathbb S)\times{\cal P}(\mathbb S)$,
$$
\sup_{f \in \Lipc}|\mu(f)-\nu(f)|=\frac{1}{2}\|\mu-\nu\|_v \, .
$$ Hence, for the 
discrete metric $\tau_c(\mathcal{M}, X)=\beta({\cal M}, \sigma(X))$
is the $\beta$-mixing coefficient between ${\cal M}$ and $\sigma(X)$ introduced in \cite{rozanov-volkonskii59beta}.
If $c$ is a distance for which $\mathbb S$ is  Polish, $\tau_c(\mathcal{M}, X)$ has 
been introduced in \cite{dedecker-prieur04coupl} and \cite{dedecker-prieur03coef}.

Applying Theorem \ref{p2} with $Q=\pr \otimes \pr_X$, we see  that this coefficient has a characteristic 
property  which is often called the {\it coupling} or {\it reconstruction} property.

\begin{cor}[reconstruction property]\label{c1} If $\Omega$ is rich enough (see Theorem \ref{p2}), there exists
a $\sigma(U)\vee \sigma(X) \vee {\cal M}$-measurable random variable  $X^*$, independent of
$\mathcal{M}$ and distributed as $X$, such that  
\begin{equation}\label{coupl}
    \tau_c(\mathcal{M}, X)=\expect{c(X, X^*)} \, .
\end{equation}
\end{cor}
If $c(x,y)= \indicatrix_{x \neq y}$, (\ref{coupl}) is given in   \cite[Corollary 4.2.5]{Berbee79book} (note that 
in Berbee's corollary,  
${\mathbb S}$ is  assumed to be standard Borel. For other proofs of Berbee's coupling, see
\cite{bryc82coupl} and \cite[Section 5.3]{rio00book}).
If $c$ is a distance for which $\mathbb S$ is a Polish space,
(\ref{coupl}) has been proved by \cite{dedecker-prieur04coupl}.

Coupling is a very useful property in the area of limit theorems and statistics.
Many authors have used Berbee's  coupling to prove various limit theorems (see for instance 
the review paper \cite{merlevede-peligrad02coupl} and the references therein)
as well as exponential inequalities (see for instance the paper \cite{doukhan-massart-rio95emp}
for Bernstein-type inequalities and applications to empirical central limit theorems).
Unfortunately, these results apply only to $\beta$-mixing sequences, but this property is very
hard to check and many simple
processes (such as iterates of maps or many non-irreducible Markov chains) 
are not $\beta$-mixing. In many cases however, this difficulty may be 
overcome by considering another distance $c$, more adapted to the problem than 
the discrete metric (typically $c$ is a norm for which ${\mathbb S}$ 
is a separable Banach space). 
The case ${\mathbb S}={\mathbb R}$ and $c(x,y)=|x-y|$, is studied in the 
paper \cite{dedecker-prieur03tau}, where many non $\beta$-mixing examples are given. 
 In this
paper the authors used the coefficients $\tau_c$ to prove Bernstein-type inequalities and 
a strong invariance principle
for partial sums. In the paper \cite[Section 4.4]{dedecker-prieur03coef} the same authors show that 
if $T$ is an uniformly expanding map  preserving a probability $\mu$ on $[0,1]$, then 
$\tau_c(\sigma(T^n), T)=O(a^n)$ for $c(x,y)=|x-y|$ and some $a$ in $[0,1[$.

The following inequality (which can be deduced from \cite[page 174]{merlevede-peligrad02coupl}) shows clearly
 that $\beta({\cal M}, \sigma(X))$ is in some sense the more restrictive 
coefficient among all the $\tau_c({\cal M}, X)$: for any $x$ in ${\mathbb S}$, we have that
\begin{equation}\label{mp}
  \tau_c({\cal M}, X) \leq 2 \int_0^{\beta({\cal M}, \sigma(X))} Q_{c(X,x)}(u) du \, ,
\end{equation}
where $Q_{c(X,x)}$ is the generalized inverse of the  function
$t \mapsto \pr(c(X,x)>t)$. In particular, if $c$ is bounded by $M$, 
$\tau_c({\cal M}, X) \leq 2M \beta({\cal M}, \sigma(X))$.

%%%%%%%%%%%%%%%%%%%%%%%%%%%%%%%%%%%%%%%%%%%%%%%%%%%%%%%%%%%%%%%%%%%%
\section{Appendix: topological and measure-theoretical complements} \label{rap}
%%%%%%%%%%%%%%%%%%%%%%%%%%%%%%%
\paragraph{Topological spaces}%
%%%%%%%%%%%%%%%%%%%%%%%%%%%%%%%
Let us recall some definitions  (see 
 \cite{schwartz73book,gardner-pfeffer84borel-measures} 
 for complements on Radon and Suslin spaces).  
 A topological space $\esp$ is said 
 to be 
 \begin{itemize}
 \item {\em regular} if, for any $x\in\esp$ and any closed subset $F$
 of $\esp$ which 
   does not contain $x$, there
   exist two disjoint open subsets $U$ and $V$ such that $x\in U$
   and $F\subset V$, 

\item {\em completely regular}
   if, for any $x\in\esp$ and any closed subset $F$ of $\esp$ which
   does not contain $x$, there
   exists a continuous function $f :\esp\ra [0,1]$ such that
   $f(x)=0$ and $f=1$ on $F$ (equivalently, $\esp$ is 
   {\em uniformizable}, that is, the topology of $\esp$ can be
   defined by a set of semidistances), 
 
 \item {\em pre-Radon} if every finite $\tau$--additive Borel measure
 on $\esp$ is 
 inner regular with respect to the compact subsets of $\esp$ 
(a Borel measure $\mu$ on $\esp$ is {\em $\tau$--additive} if, for any
 family $(F_\alpha)_{\alpha\in A}$ of closed subsets of $\esp$ such
 that 
$\foreach{\alpha,\beta\in A} \thereis{\gamma\in A} F_\gamma\subset
 F_\alpha\cap F_\beta$, we have $\mu(\cap_{\alpha\in A}F_\alpha)=\inf_{\alpha\in A}\mu(F_\alpha)$),

 \item {\em Radon} if every finite Borel measure on $\esp$ is
 inner regular with respect to the compact subsets of $\esp$,

 \item {\em Suslin}, or {\em analytic}, if there exists a continuous
   mapping from some Polish 
   space onto $\esp$, 

 \item {\em Lusin} if there exists a continuous injective mapping from
   some Polish 
   space onto $\esp$. Equivalently, $\esp$ is Lusin if there exists a
   Polish topology on $\esp$ which is finer than the given % original 
   topology of $\esp$. 
 \end{itemize}
Obviously, every Lusin space is Suslin and every Radon space is
pre-Radon. Much less obviously, every Suslin space is Radon. 
Every regular Suslin space is completely regular. 
 
Many usual spaces of Analysis are Lusin: besides all separable Banach
spaces (e.g.~$\L^p$ ($1\leq p<+\infty$), or the Sobolev spaces
$\text{\rm W}^{s,p}(\Omega)$ ($0<s<1$ and $1\leq p<+\infty$)), the spaces of distributions 
$\mathcal E'$, $\mathcal S'$, $\mathcal D'$, the space ${\mathcal
  H}(\C)$ of holomorphic functions, or the topological dual of a
Banach space, endowed with its weak$^*$--topology are Lusin. 
See \cite[pages 112--117]{schwartz73book} for many more examples.   

%%%%%%%%%%%%%%%%%%%%%%%%%%%%%%%%%%%%%%%%%%%%
\paragraph{Standard Borel spaces}
A measurable space $(\espmesgen,\tribugen)$ is said 
to be {\em standard Borel} if it is Borel-isomorphic
with some Polish space $\pol$, that is, there exists a mapping 
$f :\,\pol\ra \espmesgen$ which is one-one and onto, such that $f$ and $f^{-1}$
are measurable for $\bor{\pol}$ and $\tribugen$. 
We say that a topological space $\esp$ is {\em standard Borel} if
$(\esp,\bor{\esp})$ is standard Borel.

If $\tau_1$ and $\tau_2$ are two comparable Suslin topologies on
$\esp$, they share the same Borel sets. 
In particular, every Lusin space is standard Borel.

A useful property of standard Borel spaces is that every standard
space $\esp$ is Borel-isomorphic
with a Borel subset of $[0,1]$. This a consequence of e.g.~%~\cite[]{parthasarathy67book} or 
\cite[Theorem 15.6 and Corollary 6.5]{kechris95book}, see also
\cite{skorohod76representation} 
or \cite[Th\'eor\`eme III.20]{dellacherie-meyer75book}. 
(Actually, we have more: every
standard Borel space is countable or Borel-isomorphic with $[0,1]$. 
Thus, for standard Borel spaces, the Continuum Hypothesis holds true!)

Another useful property of standard Borel spaces is 
that, if $\esp$ is a standard Borel space, 
%a regular conditional probability always exists, that is, 
if $X :\,\esprob\mt\esp$ is a measurable mapping, and 
if $\sstribu$ is a sub-$\sigma$-algebra of $\tribu$, there exists a
regular conditional distribution $\pr_{X|\sstribu}$ 
(see e.g.~\cite[Theorem 10.2.2]{dudley02book} for the Polish case,
which immediately extends to standard Borel spaces from their definition). 
Note that, if ${\mathbb S}$ is radon, then
the distribution $\pr_X$ of $X$ is tight, that is, for every integer $n \geq 1$, there exists
a compact subset $K_n$ of ${\mathbb S}$ such that $\pr_X({\mathbb S} \setminus K_n) \geq 1/n$. Hence
one can assume without loss of generality that $X$ takes its values in $\cup_{n \geq 1} K_n$. 
If moreover ${\mathbb S}$ has metrizable compact subsets, then
$\cup_{n \geq 1} K_n$ is Lusin (and hence standard Borel), and there exists a 
regular conditional distribution $\pr_{X|\sstribu}$.
Thus, if $\esp$ is Radon with metrizable compact subsets, every element $\mu$ of $\youngs$
has an $\tribu$-measurable disintegration. Indeed, denoting
$\tribu'=\tribu\otimes\{\emptyset,\esp\}$, one only needs to consider
the conditional distribution $\pr_{X|\tribu'}$ of the random variable
$X :\,(\omega,x)\mt x$ defined on the probability space
$(\esprob\times\esp,\tribu\otimes\bor{\esp},\mu)$. 
%Actually, this result holds more generally if $\bor{\esp}$ is
%Borel-isomorphic with the Borel $\sigma$-algebra of a Radon space, this is
%an easy consequence of the disintegration theorem (e.g.~\cite{valadier73desi}). 

%\cite[Theorem 8.1]{parthasarathy67book}. 

For any $\sigma$--algebra $\tribugen$ on a set $\espmesgen$, the {\em
  universal completion} of $\tribugen$ 
is the $\sigma$-algebra 
$\tribugen^*=\cap_{\mu}\tribugen^*_\mu$, where $\mu$ runs over all
finite nonegative measures on $\tribugen$ and $\tribugen^*_\mu$ is
the $\mu$--completion of $\tribugen$. 
A subset of a topological space $\esp$ is said to be {\em
  universally measurable} if it belongs to $\bor{\esp}^*$. 
The following lemma can be deduced from e.g.~\cite[Exercise 10 page
14]{vaart-wellner96book} and the Borel-isomorphism theorem. %% between the
%% standard Borel space $\esp$ and a Borel subset of $[0,1]$. 

%%%%%%%%%%%%%%%%%%%%%%%%%%%%%%%%%%%
\begin{lem}\label{lem:modification}
Assume that $\esp$ is a standard Borel space. 
Let $X :\,\esprob\ra\esp$ be $\tribu^*$--measurable. 
Then there exists an $\tribu$--measurable modification $Y :\,\esprob\ra
\esp$ of $X$, that is, $Y$ is $\tribu$--measurable and satisfies $Y=X$ a.e.
\end{lem}
%% \preuve
%% Just take $Y=\condexpect{X}{\tribu}$. Indeed,  
%% we can assume without loss of generality that $\esp$ is Polish. Let
%% $d$ be a distance on $\esp$ which is compatible with the topology of
%% $\esp$. 
%% We denote by $\pr$ the unique extension of $\pr$ to $\tribu^*$. 
%% Let $A=\accol{X\not=Y}$ and assume that $\pr(A)>0$. 
%% We have $A=\cap_{n\geq 1}A_n$, where $A_n=\accol{d(X,Y)\geq 1/n}$,
%% thus necessarily there exists an $n\geq 1$ such that $\pr(A_n)>0$. 
%% Let $\dense=\{x_k\tq k\in\N\}$ be a countable dense subset of
%% $\esp$. We can find some $k\in\N$ such that $\pr(B)>0$, with
%% $$B=A_n\cap\bigl\{d(X,x_k)<\frac{1}{2n}\bigr\}>0.$$
%% Necessarily, we have $d(Y,x_k)>1/(2n)$ on $B$. 
%% Now, as we have $\pr=\pr_*$ on $\tribu^*$, 
%% there exists $B'\in\tribu$, with $B'\subset B$ and 
%% $\pr(B')>0$. 
%% But then we must have 
%% $$\int_{B'} \un{d(.,x_k)\leq 1/(2n)}(X)\,d\pr
%% \not=
%% \int_{B'} \un{d(.,x_k)\leq 1/(2n)}(Y)\,d\pr,$$
%% and this contradicts the definition of $Y$. 
%% \fin
%% %%%%

%%%%%%%%%%%%%%%%%%%%%%%%%%%%%%%%%%%%%%%%%%%%%%%
%\bibliographystyle{amsalpha}
%\bibliography{KR}

\begin{thebibliography}{CRdFV04}

\bibitem{Berbee79book}
Henry C.~P. Berbee, \emph{Random walks with stationary increments and renewal
  theory}, Mathematical Centre Tracts, vol. 112, Mathematisch Centrum,
  Amsterdam, 1979.

\bibitem{bryc82coupl}
W.~Bryc, \emph{On the approximation theorem of {I}. {B}erkes and {W}.
  {P}hilipp}, Demonstratio Mathematica \textbf{15} (1982), no.~3, 807--816.

\bibitem{cc-prf-valadier04book}
Charles Castaing, Paul Raynaud~de Fitte, and Michel Valadier, \emph{Young
  measures on topological spaces. {W}ith applications in control theory and
  probability theory}, Kluwer Academic Publishers, Dordrecht, 2004.

\bibitem{castaing-valadier77book}
Charles Castaing and Michel Valadier, \emph{Convex analysis and measurable
  multifunctions}, Lecture Notes in Math., no. 580, Springer Verlag, Berlin,
  1977.

\bibitem{dedecker-prieur04coupl}
J{\'e}r{\^o}me Dedecker and Cl{\'e}mentine Prieur, \emph{Couplage pour la
  distance minimale}, C. R. Math. Acad. Sci. Paris \textbf{338} (2004), no.~10,
  805--808.


\bibitem{dedecker-prieur03tau}
J\'er\^ome Dedecker and Cl\'ementine Prieur, \emph{Coupling for
$\tau$-dependent sequences and applications}, to appear in Journal of
Theoretical Probability (2003).

\bibitem{dedecker-prieur03coef}
J\'er\^ome Dedecker and Cl\'ementine Prieur, \emph{New dependence coefficients. 
Examples and applications to statistics}, accepted in Probability Theory and Related Fields
(2003).





\bibitem{dellacherie-meyer75book}
Claude Dellacherie and Paul~Andr{\'e} Meyer, \emph{Probabilit\'es et potentiel.
  {C}hapitres {I} \`a {IV}}, Hermann, Paris, 1975.


\bibitem{dobrushin70var}
    R. L. Dobru{\v{s}}in,
    \emph{Prescribing a system of random variables by 
              conditional distributions}, Theor. Probability Appl. \textbf{15} (1970), 458-486.


\bibitem{doukhan-massart-rio95emp}
Paul Doukhan, Pascal Massart, and Emmanuel Rio, \emph{Invariance principles for
  absolutely regular empirical processes}, Annales Inst. H. Poincar\'e Probab.
  Statist. \textbf{31} (1995), 393--427.



\bibitem{dudley02book}
R.~M. Dudley, \emph{Real analysis and probability}, Cambridge University Press,
  Cambridge, 2002.

\bibitem{gardner-pfeffer84borel-measures}
R.~J. Gardner and W.~F. Pfeffer, \emph{Borel measures}, Handbook of
  set-theoretic topology, North-Holland, Amsterdam, 1984, pp.~961--1043.

\bibitem{kechris95book}
Alexander~S. Kechris, \emph{Classical descriptive set theory}, Graduate Texts
  in Mathematics, no. 156, Springer Verlag, New York, 1995.

\bibitem{levin84topological}
V.~L. Levin, \emph{The problem of mass transfer in a topological space and
  probability measures with given marginal measures on the product of two
  spaces}, Dokl. Akad. Nauk SSSR \textbf{276} (1984), no.~5, 1059--1064,
  English translation: Soviet Math. Dokl. 29 (1984), no. 3, 638--643.

\bibitem{merlevede-peligrad02coupl}
Florence Merlev\`ede and Magda Peligrad, \emph{On the coupling of dependent
  random variables and applications}, Empirical process techniques for
  dependent data, Birkh{\"{a}}user, 2002, pp.~171--193.


\bibitem{rachev-ruschendorf98bookI}
S.~T. Rachev and L.~R{\" u}schendorf, \emph{Mass transportation problems.
  {V}olume {I}: Theory}, Probability and its Applications, Springer Verlag, New
  York, Berlin, 1998.



\bibitem{rio00book}
Emmanuel Rio, \emph{Th\'eorie asymptotique des processus al\'eatoires
  faiblement d\'ependants}, Math\'ematiques et Applications, no.~31, Springer,
  Berlin, Heidelberg, 2000.


\bibitem{rozanov-volkonskii59beta}
Y.~A. Rozanov and V.~A. Volkonskii, \emph{Some limit theorems for random
  functions {I}}, Teor. Verojatnost. i Primenen. \textbf{4} (1959), 186--207.

\bibitem{schwartz73book}
Laurent Schwartz, \emph{{R}adon measures on arbitrary topological spaces and
  cylindrical measures}, Tata Institute of Fundamental Research Studies in
  Mathematics, Oxford University Press, London, 1973.

\bibitem{skorohod76representation}
A.~V. Skorohod, \emph{On a representation of random variables}, Teor.
  Verojatnost. i Primenen. \textbf{21} (1976), no.~3, 645--648, English
  translation: Theor. Probability Appl. 21 (1976), no. 3, 628--632 (1977).

\bibitem{valadier73desi}
Michel Valadier, \emph{D\'esint\'egration d'une mesure sur un produit}, C. R.
  Acad. Sci. Paris S\'er. I \textbf{276} (1973), A33--A35.

\bibitem{vaart-wellner96book}
Aad~W. van~der Vaart and Jon~A. Wellner, \emph{Weak convergence and empirical
  processes. {W}ith applications to statistics}, Springer Series in Statistics,
  Springer Verlag, Berlin, 1996.

\bibitem{wheeler83survey}
Robert F. Wheeler, \emph{A survey of Baire measures and strict topologies}, 
Exposition. Math. {\bf 1} (1983), no 2, 97-190.
\end{thebibliography}
%\end{document}
\providecommand{\bysame}{\leavevmode\hbox to3em{\hrulefill}\thinspace}
\providecommand{\MR}{\relax\ifhmode\unskip\space\fi MR }
% \MRhref is called by the amsart/book/proc definition of \MR.
\providecommand{\MRhref}[2]{%
  \href{http://www.ams.org/mathscinet-getitem?mr=#1}{#2}
}
\providecommand{\href}[2]{#2}

%%%%%%%%%%%%%%%%%%%%%%%%%%%%%%%%%%%%%%%%%%%%%%%%%%%%%%%%%%%%%%%%%%%%%%%%%%%%%
\end{document}